\documentclass[11pt,letterpaper]{amsart}
\usepackage{amsmath, amsthm, amsopn, amssymb, microtype}
\usepackage{mathrsfs} 
\usepackage{mathscinet}
\usepackage{bbm}
\usepackage{enumerate, tikz, etoolbox, intcalc, geometry, caption, subcaption, afterpage}
\usepackage[english]{babel}
\usepackage{enumitem}
\usepackage[toc,page]{appendix}

\usepackage[colorlinks=true,urlcolor=blue,pdfborder={0 0 0}]{hyperref}
\hypersetup{linkcolor=[rgb]{0,0,0.6}}
\hypersetup{citecolor=[rgb]{0,0.6,0}}
\usepackage[nameinlink]{cleveref}
\crefname{theorem}{Theorem}{Theorems}
\crefname{thm}{Theorem}{Theorems}
\crefname{mainthm}{Theorem}{Theorems}
\crefname{conj}{Conjecture}{Theorems}
\crefname{lemma}{Lemma}{Lemmas}
\crefname{lem}{Lemma}{Lemmas}
\crefname{remark}{Remark}{Remarks}
\crefname{prop}{Proposition}{Propositions}
\crefname{defn}{Definition}{Definitions}
\crefname{corollary}{Corollary}{Corollaries}
\crefname{cor}{Corollary}{Corollaries}
\crefname{section}{Section}{Sections}
\crefname{figure}{Figure}{Figures}
\crefname{quest}{Question}{Questions}

\newcommand{\N}{\mathbb{N}}
\newcommand{\Z}{\mathbb{Z}}
\newcommand{\Q}{\mathbb{Q}}
\newcommand{\R}{\mathbb{R}}

\newcommand{\C}{\mathbb{C}}

\newcommand{\T}{\mathbb{T}}

\newcommand{\ind}{\mathbf{1}}

\newcommand{\K}{\mathcal{K}}
\newcommand{\Ad}{\operatorname{Ad}}
\graphicspath{ {./Pictures/} }

\newcommand{\ignore}[1]{}
\newcommand{\noopsort}[1]{}

\newtheorem{thm}{Theorem}[section]

\newtheorem{prop}[thm]{Proposition}
\newtheorem{cor}[thm]{Corollary}

\newtheorem{quest}[thm]{Question}

\theoremstyle{definition}
\newtheorem{definition}[thm]{Definition}

\newtheorem*{thm*}{Theorem}
\newtheorem*{prop*}{Proposition}

\newif\ifdraft\drafttrue

\begin{document}

\title[Criterion for Glasner property and IC-Rigidity for actions on homogeneous space]{Criterion for Glasner property and IC-Rigidity for actions on homogeneous space}

\author{Alexander Fish}
\address{School of Mathematics and Statistics, University of Sydney, NSW 2006, Australia}
\email{alexander.fish@sydney.edu.au}

\author{Shrey Sanadhya}
\address{School of Mathematics and Statistics, University of Sydney, NSW 2006, Australia}
\email{shrey.sanadhya@sydney.edu.au}


\begin{abstract} 
Let $G$ be a semi-group acting continuously on a compact metric space $X$. The action $G\curvearrowright X$ is called Glasner if for every infinite set $A\subseteq X$ and every $\varepsilon>0$, there exists $g\in G$ such that $gA$ is $\varepsilon$-dense in $X$. Let $\mathcal{K}(X)$ denote the space of nonempty compact subsets of $X$. A $G$-invariant Borel probability measure on $\mathcal{K}(X)$ is called an invariant random compact (IRC), and the action is IC-rigid if $\delta_X$ is the only IRC supported on infinite compact subsets of $X$. In this paper we give a pair-correlation criterion implying both IC-rigidity and the Glasner property. As an application, we prove both these properties for finitely generated $\Ad$-Zariski-dense semi-groups acting on compact homogeneous spaces $G/\Lambda$, where $G$ is connected and noncompact with simple Lie algebra and $\Lambda<G$ is a cocompact lattice; in particular, this applies to finitely generated Zariski-dense semi-groups of $\mathrm{SL}_d(\R)$ and $\mathrm{PSL}_d(\R)$. We also derive a Fourier-analytic criterion for semi-group actions on $\T^d$ and use it to prove IC-rigidity for actions arising from suitable matrix-valued polynomials.
\end{abstract}

\maketitle

\section{Introduction}\label{sec:Intro}

Let $G$ be a semi-group acting continuously on a compact space $X$. The action $G\curvearrowright X$ is called \textit{Glasner} if for every infinite set $Y\subset X$ and every $\varepsilon>0$ there exists $g\in G$ such that $gY$ is $\varepsilon$-dense in $X$. The terminology originates from a theorem of Glasner \cite{Glasner_1979}, which states that the multiplicative action of $\N$ on $\T=\R/\Z$ has this property. This result has since been extended in several directions. Berend and Peres \cite{Berend_Peres_1993} proved that for every infinite $A\subseteq\T$, the dilates $nA$ converge to $\T$ along a density-one set of integers. Dong \cite{Dong_2019_1,Dong_2019_2} established analogous results for higher-dimensional linear actions and actions on compact homogeneous spaces, and more recently Bulinski and Fish \cite{Bulinski_Fish_2023} obtained a quantitative version of Glasner results for irreducible linear actions with Zariski-connected Zariski closure.

A closely related notion of \textit{IC-rigidity} was introduced by Kra and Schmieding in \cite{Schmieding_Kra_2026}. Both properties concern the behavior of infinite subsets under the action: the Glasner property requires that every infinite set has arbitrarily dense images, while IC-rigidity requires that the only invariant measure supported on infinite compact subsets is $\delta_X$. We now recall the necessary background and notation to formulate the notion of IC-rigidity. 

Let $\K(X)$ denote the space of nonempty compact subsets of $X$ with Hausdorff metric. We assume that the action is surjective and finite-to-one. The action of $G$ on $X$ induces an action of $G$ on $\K(X)$. In \cite{Schmieding_Kra_2026} Kra and Schmieding introduced the notion of \emph{invariant random compacts}. A Borel probability measure $\mu$ on $\K(X)$ which is invariant under the action of $G$ on $\K(X)$ is called an invariant random compact or an IRC (see \cite[Definition 3.1]{Schmieding_Kra_2026}). An IRC $\mu$ for which
$\mu$-almost every compact subset of $X$ is infinite is called \emph{nonfinitary}.

Assuming $X$ is infinite, every infinite invariant compact set $Y \subset X$ yields a nonfinitary IRC $\delta_Y$, in particular there always exists a nonfinitary IRC, $\delta_X$. Kra and Schmieding call an action of $G$ on $X$ to be \emph{IC-rigid} if the only nonfinitary IRC is $\delta_X$ (see \cite[Definition 3.7]{Schmieding_Kra_2026}). In this paper we provide a pair-correlation condition for the action $G\curvearrowright X$ that implies IC-rigidity and the Glasner property. Let $(\mu_N)_{N\geq1}$ be a sequence of probability measures on $G$. For $f\in C(X;\C)$ define the pair-correlation kernel as follows.
\begin{equation}\label{eq:pair-kernel}
 K_{N,f}(x,y)
 :=\int_G f(gx)\overline{f(gy)}\,d\mu_N(g),
 \qquad x,y\in X.
\end{equation}

Our main result is the following.

\begin{thm}\label{thm:main}
Let $m$ be a fully supported Borel probability measure on $X$. Assume that, for every $f\in C(X;\C)$ satisfying $\int_X f\,dm=0$, every $\eta>0$, and every $x\in X$, the set
\[
\left\{ y\in X\setminus\{x\}: \limsup_{N\to\infty}|K_{N,f}(x,y)|>\eta \right\}
\]
is finite. Then:
\begin{enumerate}[label=\textup{(\alph*)},leftmargin=2em]
\item The action $G\curvearrowright X$ has the Glasner property: for every infinite set $A\subseteq X$ and every $\varepsilon>0$, there exists $g\in G$ such that $gA$ is $\varepsilon$-dense in $X$.
\item The action $G\curvearrowright X$ is IC-rigid. Equivalently, $\delta_X$ is its only nonfinitary invariant random compact.
\end{enumerate}
\end{thm}

In the context of the Glasner property for actions on compact homogeneous spaces, the best known result is due to Dong \cite{Dong_2019_2}. Dong established the Glasner property for actions of epimorphic subgroups, in particular minimal parabolic subgroups, on compact homogeneous spaces. As an application of \Cref{thm:main}, we prove both the Glasner property and IC-rigidity for finitely generated $\Ad$-Zariski-dense semi-groups acting on compact homogeneous spaces.

\begin{thm}\label{thm:simple}
Let $G$ be a connected noncompact Lie group with simple Lie algebra, and let $\Lambda<G$ be a cocompact lattice. Let $\mu$ be a finitely supported probability measure on $G$ whose support generates an $\Ad$-Zariski-dense semi-group $\Gamma\subseteq G$. Then the action $\Gamma\curvearrowright G/\Lambda$ has the Glasner property and is IC-rigid.
\end{thm}

 As a direct corollary of \Cref{thm:simple},  we show that finitely generated Zariski-dense semi-groups in $\mathrm{SL}_d(\R)$ or $\mathrm{PSL}_d(\R)$ act on their compact homogeneous spaces with the Glasner property and IC-rigidity.

\begin{cor}\label{cor:homogeneous}
Let $d\geq 2$, let
\[
G=\mathrm{SL}_d(\R) \qquad\text{or}\qquad G=\mathrm{PSL}_d(\R),
\]
and let $\Lambda<G$ be a cocompact lattice. Let $\mu$ be a finitely supported probability measure whose support generates a Zariski-dense semi-group $\Gamma<G$. Then the action $\Gamma\curvearrowright G/\Lambda$ has the Glasner property and is IC-rigid.
\end{cor}

In a recent preprint \cite{kra_Schmieding_2026_b}, Kra and Schmieding proved IC-rigidity for the action of the multiplicative semi-group $\{n^r:n\in\N\}$, for $r\geq1$, on the torus $\T$. They also proved IC-rigidity for the action on $\T^d$ of a subgroup $\Gamma<\mathrm{SL}_d(\Z)$ whose Zariski closure is Zariski-connected and whose representation on $\R^d$ is irreducible (see \cite[Theorem 1.2, Theorem 1.4]{kra_Schmieding_2026_b}). Building upon \Cref{thm:main}, we provide a Fourier-analytic criterion (see \Cref{thm:Fourier} below) for semi-group actions on the $d$-dimensional torus, $d\in\N$, which implies both IC-rigidity and the Glasner property. Our Fourier-analytic criterion also recovers the IC-rigidity results of Kra and Schmieding \cite[Theorem 1.2, Theorem 1.4]{kra_Schmieding_2026_b}, through an independently developed but closely related approach.

Recall that by a semi-group action on a torus we mean the action of a semi-group $G$ on $\T^d = \R^d/\Z^d$ by continuous group endomorphism. As before let $(\mu_N)_{N \in \N}$ be a sequence of probability measures on $G$. The averaged orbit probability measure of $z\in \T^d$ (corresponding to $\mu_N$), is given by 
\[
\sigma_{N,z}
:=
\int_G \delta_{gz}\, d\mu_N(g).
\]
For $t \in \R$, we set $e(t):=e^{2\pi i t}$. For $m \in \Z^d$, the $m$-th Fourier coefficient of the probability measure $\sigma_{N,z}$ is given by
\[
\widehat{\sigma_{N,z}}(m)
=
\int_{\mathbb T^d}
e(\langle m,x\rangle)\,d\sigma_{N,z}(x)
=
\int_G e(\langle m,gz\rangle)\,d\mu_N(g).
\]
We show that \Cref{thm:main}, together with a density argument, gives the following.

\begin{thm}\label{thm:Fourier}
Assume that there exists a sequence of probability measures $(\mu_N)_{N\geq1}$ on $G$ such that, for every $m\in\Z^d\setminus\{0\}$ and every $\delta>0$, the set  
\begin{equation}\label{eq:criteria}
\left\{z\in\T^d: \limsup_{N\to\infty} \left| \int_G e\bigl(\langle m,gz\rangle\bigr)\,d\mu_N(g) \right| \geq\delta \right\}
\end{equation} is finite. Then:
\begin{enumerate}[label=\textup{(\alph*)},leftmargin=2em]
\item The action $G\curvearrowright\T^d$ has the Glasner property: for every infinite set $A\subseteq\T^d$ and every $\varepsilon>0$, there exists $g\in G$ such that $gA$ is $\varepsilon$-dense in $\T^d$. 
\item The action $G\curvearrowright\T^d$ is IC-rigid. Equivalently, $\delta_{\T^d}$ is its only nonfinitary invariant random compact.
\end{enumerate} 
\end{thm}

In \cite{Bulinski_Fish_2023_a} Bulinski and Fish proved Glasner property for actions arising from matrix-valued polynomials. They showed that for $A(x)\in M_d(\Z[x])$ the condition
\[
v\cdot(A(x)-A(0))w\neq0
\]
for every $v,w\in\Z^d\setminus\{0\}$ is equivalent to the set $\{A(n):n\in\Z\}$ having the uniform Glasner property on $\T^d$. As an application of \Cref{thm:Fourier} (b), we prove IC-rigidity for the multiplicative semi-group generated by such a matrix-valued polynomial. In private communication, Bryna Kra and Scott Schmieding informed us that this result can also be obtained using their methods in \cite{kra_Schmieding_2026_b}.

\begin{thm}\label{thm:polynomial-matrix}
Let $A(x)\in M_d(\Z[x])$ satisfy 
\[
v\cdot (A(x)-A(0))w\neq 0
\]
for every $v,w\in\Z^d\setminus\{0\}$, and suppose that $\det A(x)\not\equiv 0$. Let $G_A$ be the multiplicative semi-group generated by
\[
\{A(n):n\in\Z,\ \det A(n)\neq 0\}.
\]
The action $G_A\curvearrowright\T^d$ by endomorphisms, given by
\[
g\cdot z=gz\pmod{\Z^d},
\] is IC-rigid. Equivalently, $\delta_{\T^d}$ is the only nonfinitary invariant random compact for this action.
\end{thm}

\subsection{Open questions}

We conclude with a few open questions. For relatively prime integers $p,q\geq2$, let
\[
S_{p,q}:=\{p^aq^b:a,b\geq0\}
\]
be the multiplicative semi-group generated by $p$ and $q$, acting on $\T$ by dilations. It is natural to ask whether this action has the Glasner property. A closely related question concerning IC-rigidity was posed by Kra and Schmieding in \cite{Schmieding_Kra_2026}. In a recent preprint \cite{Schmieding_2026}, Schmieding proved that if $\mu$ is an IRC for this action and $\mu$-almost every compact set $Y$ has positive Hausdorff dimension, then $\mu$-almost every $Y$ is equal to $\T$. In other words, $\mu=\delta_{\T}$.

\begin{quest}
Let $p,q\geq2$ be relatively prime. Does the action $S_{p,q}\curvearrowright\T$ have the Glasner property? Is this action IC-rigid?
\end{quest}

Another question posed by Kra and Schmieding in \cite{Schmieding_Kra_2026}, asks whether there exists an IC-rigid action of an abelian group, and in particular an IC-rigid $\Z$-action. Recently, Kanigowski, Kasprzak and Lorenzo-Laguno \cite{Kanigowski_Kasprzak_LorenzoLaguno_2026} answered this question positively by proving that mixing rank-one $\Z$-actions are IC-rigid. The IC-rigid $\Z$-actions constructed there are realized as subshifts, and hence the underlying spaces are totally disconnected. On the other hand, Dong \cite[Theorem 3.1]{Dong_2019_1} proved that no $\Z$-action on a locally connected compact metric space has the Glasner property. These results motivate the following questions.

\begin{quest}
\leavevmode
\begin{enumerate}
    \item[(a)] Does there exist, for some $d\geq 2$, a $\Z^d$-action on an infinite connected compact metric space with the Glasner property?
    
    \item[(b)] Does there exist, for some $d\geq 1$, an IC-rigid $\Z^d$-action on an infinite connected compact metric space?
\end{enumerate}
\end{quest}

\subsection{Structure of the paper}
In \Cref{sec:Back}, we provide the basic background and definitions. In \Cref{sec:Main}, using a fundamental result concerning the density of measures, we prove the pair-correlation criterion, \Cref{thm:main}. In \Cref{sec:Z_dense}, we recall an equidistribution result for homogeneous dynamics due to Benoist--Quint \cite{Benoist_Quint_2013} (see \Cref{thm:BQ-homogeneous}). Using this equidistribution result together with \Cref{thm:main}, we prove \Cref{thm:simple}. In the same section, as a corollary we establish the Glasner property and IC-rigidity for Zariski-dense semi-groups of $\mathrm{SL}_d(\R)$ and $\mathrm{PSL}_d(\R)$ acting on compact homogeneous spaces, as stated in \Cref{cor:homogeneous}. In \Cref{sec:Torus}, using a density argument together with \Cref{thm:main}, we prove the Fourier-analytic criterion for torus actions, \Cref{thm:Fourier}. Finally in the same section, we use \Cref{thm:Fourier} to prove \Cref{thm:polynomial-matrix}.

{\bf Acknowledgement.} The authors would like to thank Bryna Kra and Scott Schmieding for communicating their results and for their encouragement. Both authors were supported by the Australian Research Council through grant DP240100472.

{\bf AI statement.}  AI was used in the paper to obtain the relevant version of the equidistribution theorem for homogeneous spaces due to Benoist--Quint \cite{Benoist_Quint_2013} (see \Cref{thm:BQ-homogeneous}) in a form suitable for our purposes. Apart from this, AI was used to audit the paper, identify relevant references, and correct editorial typos.

\section{Background and notation}\label{sec:Back} Let $(X,d)$ be a compact metric space and $G$ be a countable semi-group or group. By $(X,G)$ we denote a topological dynamical system where $G$ acts on $X$ continuously. In the case where $G$ is a semi-group, we assume that the action is surjective and finite-to-one. 

\begin{definition} A continuous action of a semi-group $G$ on a compact space $X$ is called \emph{Glasner} if for every infinite set $Y \subset X$ and $\varepsilon>0$ there exists $g \in G$ such that $gY$ is $\varepsilon$-dense in $X$.
\end{definition}
Following the work of Glasner \cite{Glasner_1979}, this property has been extensively studied in many different contexts. We refer our interested readers to \cite{Berend_Peres_1993,  KellyLe2013, Dong_2019_1, Dong_2019_2}.

As mentioned before, we denote by $\K(X)$ the space of nonempty compact subsets of $X$ which is endowed with Hausdorff metric. For $K, L \in \K(X)$ set
\begin{equation}\label{eq:Haus}
d_H(K,L)
=
\max\left\{
\sup_{x\in K} d(x,L),
\sup_{y\in L} d(y,K)
\right\},
\end{equation} 
where $d(x,L):=\inf_{y\in L}d(x,y)$. The action of $G$ on $X$ induces an action of $G$ on $\mathcal{K}(X)$ given by $g\cdot K = g(K)$. Using the notation of Kra and Schmieding for $k \in \N$,  we denote by $\K_{\leq k}(X)$ the collection of all subsets of $X$ of size at most $k$ and set
\[
\K_{\mathrm{fin}}(X)
=
\bigcup_{k=1}^{\infty}\mathcal{K}_{\leq k}(X),
\qquad
\K_{\mathrm{inf}}(X)
=
\mathcal{K}(X)\setminus \mathcal{K}_{\mathrm{fin}}(X).
\]
Since $\K_{\leq k}(X)$ is a closed set for each $k \in \N$, it follows that $\K_{\mathrm{fin}}(X)$ is a $F_{\sigma}$ set, hence it is Borel. This implies that $\K_{\mathrm{inf}}(X)$ is also Borel. Observe that $\K_{\mathrm{fin}}(X)$ is $G$-invariant and the finite-to-one property of the action implies that $\K_{\mathrm{inf}}(X)$ is $G$-invariant as well. Now we formally define:

\begin{definition} (see \cite[Definition 3.1]{Schmieding_Kra_2026}) An \emph{invariant random compact (IRC)} for a topological dynamical system $(X,G)$, is a Borel probability measure $\mu$ on $\mathcal{K}(X)$ which is invariant under the induced action of $G$ on $\K(X)$.  
\end{definition}

Observe that an invariant random compact $\mu$ can be decomposed into its \emph{finitary and nonfinitary} parts i.e. $\mu=\mu_F+\mu_I$, where the measures $\mu_F,\mu_I$ on $\K(X)$ are defined as,
\[
\mu_F(\mathcal{A})=\mu\bigl(\mathcal{A}\cap \K_{\mathrm{fin}}(X)\bigr)
\qquad\text{and}\qquad
\mu_I(\mathcal{A})=\mu\bigl(\mathcal{A}\cap \K_{\mathrm{inf}}(X)\bigr)
\]
for a Borel set $\mathcal{A} \subseteq \K(X)$.

\begin{definition} (see \cite[Definition 3.7]{Schmieding_Kra_2026})
The action of $G$ on $X$ is \emph{infinite compact rigid (or IC-rigid)} if the only nonfinitary IRC is $\delta_X$. In other words for every IRC $\mu$, the decomposition of $\mu$ into its finitary and nonfinitary parts has the form
\[
\mu=\mu_F+c\cdot\delta_X
\]
for some $c\geq 0$. 
\end{definition}

\section{Proof of \Cref{thm:main}}\label{sec:Main} In this section we prove \Cref{thm:main}. To do so we first need an elementary result showing that small averages of finitely many mean-zero continuous functions force dense support.

\begin{prop}\label{prop:density-detection}
Let $m$ be a fully supported Borel probability measure on $X$. For every $\varepsilon>0$, there exist $q\in \N$, zero-mean functions $f_1,\ldots,f_q\in C(X;\C)$, and constants $c_1,\ldots,c_q>0$ such that the following holds. If $\nu$ is a Borel probability measure on $X$ satisfying
\[
\left|\int_X f_j\,d\nu\right|<c_j
\qquad\text{for every }1\leq j\leq q,
\]
then $\operatorname{supp}\nu$ is $\varepsilon$-dense in $X$.
\end{prop}

\begin{proof}
Choose $y_1,\ldots,y_q\in X$ such that the open balls
\[
U_j:=B(y_j,\varepsilon/2),\qquad j=1,\ldots,q,
\]
cover $X$. For every $j\in\{1,\ldots,q\}$, choose a nonzero non-negative bump function $\phi_j$ supported in $U_j$. Since $m$ is fully supported,
\[
a_j:=\int_X\phi_j\,dm>0.
\]
For every $j\in\{1,\ldots,q\}$, set $f_j:=\phi_j-a_j$ and $c_j:=\frac{a_j}{2}$. Using the fact that $m$ is a probability measure, for every $j\in\{1,\ldots,q\}$, we get $\int_Xf_j\,dm=0$. Let $\nu$ be a Borel probability measure on $X$ such that
\[
\left|\int_Xf_j\,d\nu\right|<c_j \qquad\text{for every }1\leq j\leq q.
\]
Then
\[
\int_X\phi_j\,d\nu
=
a_j+\int_Xf_j\,d\nu
>
\frac{a_j}{2}>0.
\]
Since $\phi_j$ vanishes outside $U_j$ and $\int_{X} \phi_j\,d\nu>0,$ it follows that $\operatorname{supp}\nu\, \cap \,U_j\neq\varnothing$ for every $j \in \{1,\ldots,q\}$. Since $U_j$ are balls of radius $\varepsilon/2$ that cover $X$, it follows that $\operatorname{supp}\nu$ is $\varepsilon$-dense in $X$. 
\end{proof}

Now we are ready to prove \Cref{thm:main}.

\begin{proof}[Proof of \Cref{thm:main}]
We first establish a pair-correlation estimate that will be used in proving both part (a) and (b) of the statement. Fix $\varepsilon>0$, and let $q\in\N$, zero-mean functions $f_1,\ldots,f_q\in C(X;\C)$, and constants $c_1,\ldots,c_q>0$ be as given by \Cref{prop:density-detection}. Let $A\subseteq X$ be infinite and let $\eta>0$. Set $C:=\sum_{j=1}^q\|f_j\|_\infty^2$ and $\delta:=\frac{\eta}{2q}$. Choose $M\in\N$ such that $M>\frac{2C}{\eta}$. For $1\leq j\leq q$ and $x\in X$, let
\[
R_{j,\delta}(x)
:=
\left\{
y\in X\setminus\{x\}:
\limsup_{N\to\infty}|K_{N,f_j}(x,y)|>\delta
\right\}.
\]
By the assumption of \Cref{thm:main}, $R_{j,\delta}(x)$ is finite for every $j$ and every $x\in X$. Note that for any $x \neq y \in X$, we have $K_{N,f_j}(y,x) = \overline{K_{N,f_j}(x,y)}$ and thus $|K_{N,f_j}(y,x)|=|K_{N,f_j}(x,y)|$. Thus the expression 
\[
\limsup_{N\to\infty}|K_{N,f_j}(x,y)|>\delta
\]
is symmetric in $x$ and $y$. Since $A$ is infinite and each set $R_{j,\delta}(x)$ is finite, we can choose distinct points $x_1,\ldots,x_M\in A$, inductively so that, for every distinct $i,\ell\in\{1,\ldots,M\}$ and every $1\leq j\leq q$,
\[
\limsup_{N\to\infty} |K_{N,f_j}(x_i,x_\ell)| \leq\delta.
\]
We define
\[
\tau_M:=\frac1M\sum_{i=1}^M\delta_{x_i}.
\]
Note that $\tau_M$ is the probability measure on $X$ obtained by averaging the Dirac measures at $x_1,\ldots,x_M$. For $g\in G$, the pushforward of $\tau_M$ is given by
\[
g_*\tau_M = \frac1M\sum_{i=1}^M\delta_{gx_i}.
\]
For $1\leq j\leq q$, set
\[
F_j(g) := \int_Xf_j\,d(g_*\tau_M) = \frac1M\sum_{i=1}^Mf_j(gx_i).
\]
Expanding the square and integrating over $G$ with respect to $\mu_N$ gives
\[
\begin{aligned}
\int_G|F_j(g)|^2\,d\mu_N(g) &= \frac1{M^2} \sum_{i,\ell=1}^M \int_G f_j(gx_i)\overline{f_j(gx_\ell)} \,d\mu_N(g)\\ &= \frac1{M^2} \sum_{i,\ell=1}^M K_{N,f_j}(x_i,x_\ell).
\end{aligned}
\]
For the diagonal terms $i=\ell$,
\[
|K_{N,f_j}(x_i,x_i)| = \int_G|f_j(gx_i)|^2\,d\mu_N(g) \leq\|f_j\|_\infty^2.
\]
Since there are $M$ diagonal terms, using triangle inequality we get
\[
\int_G|F_j(g)|^2\,d\mu_N(g) \leq \frac{\|f_j\|_\infty^2}{M} + \frac1{M^2} \sum_{i\neq\ell} |K_{N,f_j}(x_i,x_\ell)|.
\]
By the choice of $x_1,\ldots,x_M$, it follows that
\[
\limsup_{N\to\infty}
\int_G|F_j(g)|^2\,d\mu_N(g) \leq \frac{\|f_j\|_\infty^2}{M} + \frac{M(M-1)}{M^2}\delta \leq \frac{\|f_j\|_\infty^2}{M} +\delta.
\]
Summing over $1\leq j\leq q$ yields
\begin{equation}\label{eq:common-pair-estimate}
\limsup_{N\to\infty}\int_G\sum_{j=1}^q|F_j(g)|^2\,d\mu_N(g) \leq \sum_{j=1}^q\limsup_{N\to\infty}\int_G|F_j(g)|^2\,d\mu_N(g) \leq \frac{C}{M}+q\delta<\eta.
\end{equation}
Thus, for every infinite $A\subseteq X$, every finite collection of zero-mean functions $f_1,\ldots,f_q$, and every $\eta>0$, there is a finite set of distinct points $x_1,\ldots,x_M\in A$ for which \eqref{eq:common-pair-estimate} holds.

We now show that the action $G\curvearrowright X$ has the Glasner property. Fix an infinite set $A\subseteq X$ and $\varepsilon>0$. Let $q\in\N$, zero-mean functions $f_1,\ldots,f_q\in C(X;\C)$, and constants $c_1,\ldots,c_q>0$ be as given by \Cref{prop:density-detection}, and set
\[
 c:=\min_{1\leq j\leq q}c_j>0.
\]
For $\eta=c^2/2$, choose $x_1,\ldots,x_M\in A$ such that, for the corresponding $\tau_M$ and $F_j$, \eqref{eq:common-pair-estimate} holds. Thus, for a sufficiently large $N$,
\[
\int_G\sum_{j=1}^q|F_j(g)|^2\,d\mu_N(g)<c^2.
\]
Since $\mu_N$ is a probability measure, there exists $g\in G$ such that $\sum_{j=1}^q|F_j(g)|^2<c^2$ otherwise the integral would be at least $c^2$. Hence
\[
\left| \int_Xf_j\,d(g_*\tau_M) \right| = |F_j(g)|<c\leq c_j \qquad\text{for every }1\leq j\leq q.
\]
Thus \Cref{prop:density-detection} implies that
\[
\operatorname{supp}(g_*\tau_M) = \{gx_1,\ldots,gx_M\} \subseteq gA
\]
is $\varepsilon$-dense in $X$. Hence $gA$ is also $\varepsilon$-dense in $X$. This proves the Glasner property.

We now show that the action $G\curvearrowright X$ is IC-rigid. Let $\nu$ be a nonfinitary invariant random compact. For $\varepsilon>0$, set
\[
\K_\varepsilon := \{K\in\K(X):d_H(K,X)<\varepsilon\}.
\]
Note that the map $K\mapsto d_H(K,X)$ is continuous hence $\K_\varepsilon$ is an open subset of $\K(X)$. We show that for every infinite $K\in\K(X)$,
\begin{equation}\label{eq:abstract-dense-images}
\lim_{N\to\infty}
\mu_N\{g\in G:gK\in\K_\varepsilon\} = 1.
\end{equation}
Fix such a $K$ and let $q\in\N$, zero-mean functions $f_1,\ldots,f_q\in C(X;\C)$, and constants $c_1,\ldots,c_q>0$ be as given by \Cref{prop:density-detection}. We set $c:=\min_{1\leq j\leq q}c_j>0$. For any $\eta>0$, choose $x_1,\ldots,x_M\in K$ such that, for the corresponding $\tau_M$ and $F_j$, \eqref{eq:common-pair-estimate} holds. Thus,
\begin{equation}\label{eq:abstract-bound-eta}
\limsup_{N\to\infty} \int_G\sum_{j=1}^q|F_j(g)|^2\,d\mu_N(g) < \eta.
\end{equation}
As we noted earlier, if $\sum_{j=1}^q|F_j(g)|^2<c^2$, then
\[
\left|\int_Xf_j\,d(g_*\tau_M)\right| = |F_j(g)| < c \leq c_j \qquad\text{for every }1\leq j\leq q.
\]
Hence, by \Cref{prop:density-detection}, $\operatorname{supp}(g_*\tau_M) = \{gx_1,\ldots,gx_M\} \subseteq gK$ is $\varepsilon$-dense in $X$. In other words, $gK\in\K_\varepsilon$. Thus,
\[
\{g\in G:gK\notin\K_\varepsilon\} \subseteq \left\{ g\in G: \sum_{j=1}^q|F_j(g)|^2\geq c^2 \right\}.
\]
Taking measure $\mu_N$ of both sides and using Markov's inequality we conclude
\[
\mu_N\{g\in G:gK\notin\K_\varepsilon\} \leq \frac1{c^2} \int_G \sum_{j=1}^q|F_j(g)|^2\,d\mu_N(g).
\]
Taking $\limsup_{N\to\infty}$ of both sides and using \eqref{eq:abstract-bound-eta}, we get
\[
\limsup_{N\to\infty} \mu_N\{g\in G:gK\notin\K_\varepsilon\} \leq \frac{\eta}{c^2}.
\]
Since $\eta>0$ is arbitrary, letting $\eta\rightarrow0$ we obtain
\[
\limsup_{N\to\infty} \mu_N\{g\in G:gK\notin\K_\varepsilon\} = 0,
\]
which proves \eqref{eq:abstract-dense-images}. For $K\in\K(X)$, $N\in\N$ and $\varepsilon>0$, set
\[
H_{N,\varepsilon}(K) := \mu_N\{g\in G:gK\in\K_\varepsilon\} = \int_G\ind_{\K_\varepsilon}(gK)\,d\mu_N(g).
\]
Since $\mu_N$ is a probability measure, $0\leq H_{N,\varepsilon}(K)\leq1$ for $K\in\K(X)$. Also, using the fact that $\nu$ is nonfinitary, \eqref{eq:abstract-dense-images} implies that $H_{N,\varepsilon}(K)\rightarrow 1$ for $\nu$-almost every $K\in\K(X)$. Hence, by the dominated convergence theorem we get,
\begin{equation}\label{eq:abstract-DCT-density}
\lim_{N\to\infty} \int_{\K(X)} H_{N,\varepsilon}(K)\,d\nu(K) = 1.
\end{equation}
Using Fubini's theorem and the $G$-invariance of $\nu$, we get
\[
\begin{aligned} 
\int_{\K(X)} H_{N,\varepsilon}(K)\,d\nu(K)
&=
\int_G \int_{\K(X)} \ind_{\K_\varepsilon}(gK)\,d\nu(K)\,d\mu_N(g)\\
&=
\int_G \nu\{K\in\K(X):gK\in\K_\varepsilon\}\,d\mu_N(g)\\
&=
\int_G\nu(\K_\varepsilon)\,d\mu_N(g) = \nu(\K_\varepsilon).
\end{aligned}
\]
Letting $N\rightarrow\infty$ on both sides and using \eqref{eq:abstract-DCT-density}, we get $\nu(\K_\varepsilon)=1$, for every $\varepsilon>0$. In particular,
\[
\nu\left( \bigcap_{n=1}^\infty\K_{1/n} \right) = 1.
\]
But
\[
\bigcap_{n=1}^\infty\K_{1/n} = \{K\in\K(X):d_H(K,X)=0\} = \{X\}.
\]
Hence $\nu(\{X\})=1$, thus $\nu=\delta_X$. Thus the action $G\curvearrowright X$ is IC-rigid.

\end{proof}

\section{Zariski-dense actions on compact homogeneous spaces}\label{sec:Z_dense}

In this section, we prove \Cref{thm:simple}. Let $X=G/\Lambda$ and let $m$ denote the $G$-invariant probability measure on $X$. For a semi-group $\Gamma\subseteq G$, we denote by 
\[
\Gamma^\Delta:=\{(\gamma,\gamma):\gamma\in\Gamma\}\subseteq G\times G,
\] the diagonal semi-group associated with $\Gamma$. To prove \Cref{thm:simple}, we will apply an equidistribution result of Benoist--Quint \cite{Benoist_Quint_2013} to the action $\overline{\Gamma^\Delta} \curvearrowright X \times X$ where $\Gamma \subseteq G$ is as in the statement and $\overline{\Gamma^\Delta}$ denotes the closure of $\Gamma^\Delta$ in $G \times G$. Using the equidistribution result, we will show that the action $\Gamma \curvearrowright X$ satisfies the hypothesis of \Cref{thm:main}, hence it is Glasner and IC-rigid. 

\subsection{Equidistribution result of Benoist--Quint} We now recall an equidistribution result due to Benoist--Quint that will be used in the proof of \Cref{thm:simple}. Let $H$ be a real Lie group, let $L<H$ be a cocompact lattice, and let $\vartheta$ be a compactly supported Borel probability measure on $H$. Let
\[
\Sigma:=\overline{\langle\operatorname{supp}\vartheta\rangle_+} \subseteq H
\]
be the closed sub-semi-group of $H$ generated by $\operatorname{supp}\vartheta$. Thus $\Sigma$ is a compactly generated closed sub-semi-group of $H$. Assume that the Zariski closure of $\Ad(\Sigma)$ is semisimple with no compact factor. For $k\geq1$, let $\vartheta^{*k}$ denote the $k$-fold convolution of $\vartheta$, that is, the distribution of a product $h_k\cdots h_1$ of $k$ independent elements of $H$ chosen according to $\vartheta$. For completeness we set $\vartheta^{*0}:=\delta_e$, where $e$ is the identity element of $H$. For $z\in H/L$ and $k\geq0$, define
\[
\vartheta^{*k}*\delta_z := \int_H\delta_{hz}\,d\vartheta^{*k}(h).
\]
Thus $\vartheta^{*k}*\delta_z $ is a probability measure on $H/L$. It is the $k$-step random-walk distribution starting from $z$. Recall that a closed subset $Y$ of a homogeneous space $H/L$ is called a \textit{finite-volume homogeneous subspace} if its stabilizer $H_Y$ acts transitively on $Y$ and there exists an $H_Y$-invariant Borel probability measure $\nu_Y$ on $Y$. We call $\nu_Y$ the \textit{homogeneous probability measure} on $Y$. We use the following consequence of results of Benoist--Quint \cite[Theorems 1.1, 1.2 and 1.7]{Benoist_Quint_2013}, stated in the form adapted to the probability measure $\vartheta$ above.

\begin{thm} \cite{Benoist_Quint_2013} \label{thm:BQ-homogeneous}
Let $H$, $L$, $\vartheta$ and $\Sigma$ be as above and for $z\in H/L$, set $Y_z:=\overline{\Sigma z}$. Then $Y_z$ is a $\Sigma$-invariant, $\Sigma$-ergodic finite-volume homogeneous subspace. Let $\nu_z$ denote the homogeneous probability measure on $Y_z$. Then
\[
\frac1N\sum_{k=0}^{N-1}\vartheta^{*k}*\delta_z \longrightarrow \nu_z
\]
in the weak-* topology. If, in addition, the centralizer of $\Sigma$ in $H$ is discrete, then for every sequence $(z_n)_{n\in\N}$ in $H/L$, the corresponding homogeneous subspaces $Y_{z_n}:=\overline{\Sigma z_n}$ have a subsequence $(Y_{z_{n_m}})_{m\in\N}$ converging to a $\Sigma$-invariant, $\Sigma$-ergodic finite-volume homogeneous subspace $Y$. Moreover, $Y_{z_{n_m}}\subseteq Y$ for all sufficiently large $m\in\N$.
\end{thm}

\subsection{Proof of \Cref{thm:simple}} In this subsection we prove \Cref{thm:simple} using \Cref{thm:main} and \Cref{thm:BQ-homogeneous}. 

\begin{proof}[Proof of \Cref{thm:simple}] 
Let $\mu$ be a finitely supported probability measure on $G$ whose support generates an $\Ad$-Zariski-dense semi-group $\Gamma\subseteq G$. For $k\geq1$, let $\mu^{*k}$ denote the $k$-fold convolution of $\mu$, that is, the distribution of a product $g_k\cdots g_1$ of $k$ independent elements chosen according to $\mu$. For $N\geq1$, define
\[
\mu_N:=\frac1N\sum_{k=1}^{N}\mu^{*k}.
\]
Thus $\mu_N$ is a probability measure on $\Gamma$. We will show that the sequence $(\mu_N)_{N\geq1}$ satisfies the hypothesis of \Cref{thm:main}. We set
\[
H:=G\times G,\qquad L:=\Lambda\times\Lambda,
\]
and identify $H/L = X\times X =(G\times G)/(\Lambda\times\Lambda)$. Since $\Lambda$ is cocompact in $G$, the lattice $L$ is cocompact in $H$. Let $\Delta$ be the diagonal map given by
\[
\Delta:G\longrightarrow G\times G,\qquad \Delta(g):=(g,g). 
\]
Let $\vartheta:=\Delta_*\mu$ be the pushforward of the measure $\mu$ under the diagonal map. Since $\mu$ is finitely supported probability measure on $G$, $\vartheta$ is a compactly supported probability measure on $H$. Moreover, its support generates the diagonal semi-group $\Gamma^\Delta$. Hence, using the notation of \Cref{thm:BQ-homogeneous}, we have $\Sigma=\overline{\Gamma^\Delta}$. Observe that, for $k \in \N$, $\vartheta^{*k}=\Delta_*(\mu^{*k})$. In other words $\vartheta^{*k}$ is the pushforward of $\mu^{*k}$. 

Let $\mathfrak g:=\operatorname{Lie}(G)$, which is simple by assumption. Hence it follows that $\operatorname{Lie}(H)=\mathfrak g\oplus\mathfrak g$ and for $(\gamma,\gamma) \in H$, $\Ad_H(\gamma,\gamma) = (\Ad(\gamma),\Ad(\gamma))$. Since $\Sigma=\overline{\Gamma^\Delta}$, the Zariski closure of $\Ad_H(\Sigma)$ is the diagonal copy of the Zariski closure of $\Ad(\Gamma)$. Since $\Gamma$ is $\Ad$-Zariski dense in $G$, this is the diagonal copy of the Zariski closure of $\Ad(G)$. Since $\mathfrak g$ is simple and $G$ is noncompact, the Zariski closure of $\Ad_H(\Sigma)$ is semisimple with no compact factor. Hence the Zariski-closure hypothesis of \Cref{thm:BQ-homogeneous} is satisfied.

We claim that the centralizer of $\Sigma$ in $H$ is discrete. Let $(U,V)\in\operatorname{Lie}(Z_H(\Sigma))$, since $(U,V)$ belongs to the Lie algebra of the centralizer it is fixed under the adjoint action of every element of $\Sigma$, that is, $\Ad_H(\sigma)(U,V)=(U,V)$ for every $\sigma\in\Sigma$. In particular, since $\Gamma^\Delta\subseteq\Sigma$, for every $\gamma\in\Gamma$ we have
\[
\Ad_H(\gamma,\gamma)(U,V) = (\Ad(\gamma)U,\Ad(\gamma)V) = (U,V).
\]
Hence $\Ad(\gamma)U=U$ and $\Ad(\gamma)V=V$ for every $\gamma\in\Gamma$. By the $\Ad$-Zariski density of $\Gamma$, the same holds for every $g\in G$. Now, for $Z\in\mathfrak g$, taking $g=\exp(tZ)$ gives
\[
\Ad(\exp(tZ))U=U, \qquad \Ad(\exp(tZ))V=V
\]
for every $t$. Thus $\Ad(\exp(tZ))U$ is constant in $t$. Differentiating at $t=0$, and using the fact that 
\[
\left.\frac{d}{dt}\right|_{t=0}\Ad(\exp(tZ))U=[Z,U]
\]
we get $[Z,U] = 0$ for every $Z\in\mathfrak g$. Similarly we get $[Z,V] = 0$ for every $Z\in\mathfrak g$. In other words $U$ and $V$ belong to the center of $\mathfrak g$. Since the center of the simple Lie algebra $\mathfrak g$ is trivial, $U=V=0$. Therefore $\operatorname{Lie}(Z_H(\Sigma))=0$. Since $Z_H(\Sigma)$ is a closed subgroup of the Lie group $H$, this implies that $Z_H(\Sigma)$ is discrete. Hence the hypothesis of \Cref{thm:BQ-homogeneous} holds. 

Now fix $f\in C(X;\C)$ satisfying $\int_X f\,dm=0$, and for $u,v \in X$, set $\Phi_f(u,v):=f(u)\overline{f(v)}$ for convenience. For $x,y\in X$, set
\[
Y_{x,y}:=\overline{\Sigma(x,y)}
\]
and let $\nu_{x,y}$ denote its homogeneous probability measure. By \Cref{thm:BQ-homogeneous},
\begin{equation} \label{eq:w_star}
    \frac1N\sum_{k=0}^{N-1} \vartheta^{*k}*\delta_{(x,y)} \longrightarrow \nu_{x,y}
\end{equation}
in the weak-* topology. Hence we get,
\begin{equation}\label{eq:conv_K_N_f}
K_{N,f}(x,y) = \int_G f(gx)\overline{f(gy)}\,d\mu_N(g) = \int_{X\times X}\Phi_f\,d\left(\frac1N\sum_{k=1}^{N}\vartheta^{*k}*\delta_{(x,y)}\right)\xrightarrow[N\to\infty]{} \int_{Y_{x,y}}\Phi_f\,d\nu_{x,y}.    
\end{equation}
Here we used the fact that averages from $k=1$ to $N$ have the same limit as \eqref{eq:w_star}, since the two averages differ only in the $k=0$ and $k=N$ terms, each divided by $N$. Fix $x\in X$ and $\eta>0$. We claim that
\[
E_{f,\eta}(x):= \left\{ y\in X\setminus\{x\}: \limsup_{N\to\infty}|K_{N,f}(x,y)|>\eta \right\}
\]
is finite. To see this assume by contradiction that $E_{f,\eta}(x)$ is infinite. For $n \in \N$, choose distinct $y_n \in E_{f,\eta}(x)$ and set
\[
Y_n:=Y_{x,y_n} =\overline{\Sigma(x,y_n)}, \qquad
\nu_n:=\nu_{x,y_n}.
\] Here for $n \in \N$, $\nu_n$ denotes the homogeneous probability measure on $Y_n$. Then by definition of $E_{f,\eta}(x)$ and \eqref{eq:conv_K_N_f} we have
\[
\left| \int_{Y_n}\Phi_f\,d\nu_n \right|>\eta \qquad\text{for every }n.
\]
By \Cref{thm:BQ-homogeneous}, after passing to a subsequence, there exists a $\Sigma$-invariant, $\Sigma$-ergodic finite-volume homogeneous subspace $Y\subseteq X\times X$ such that $Y_n\subseteq Y$ for all sufficiently large $n$, and $\nu_n\longrightarrow\nu_Y$ in the weak-* topology. This implies 
\begin{equation}\label{eq:_Phi_contra}
\left| \int_Y\Phi_f\,d\nu_Y \right| \geq\eta.
\end{equation}
Since $(x,y_n)\in Y_n\subseteq Y$ for all sufficiently large $n$, the fiber $F:=Y\cap(\{x\}\times X)$ is infinite. Let $H_Y$ be the stabilizer of $Y$ in $H$. Since $Y$ is a finite-volume homogeneous subspace, $H_Y$ acts transitively on $Y$ and preserves the homogeneous probability measure $\nu_Y$. We set
\[
G_x:=\{g\in G:gx=x\}.
\]
Since $X=G/\Lambda$, the group $G_x$ is conjugate to $\Lambda$, and hence is discrete. The group
\[
K:=H_Y\cap(G_x\times G)
\]
acts transitively on $F$. Assume that $K$ were discrete. Since $H$ is a second-countable Lie group, $K$ would be countable. As $K$ acts transitively on $F$, the set $F$ would be countable and homogeneous under the action of $K$. Moreover, $F=Y\cap(\{x\}\times X)$ is compact, hence $F$ would be a countable compact homogeneous space. Recall that, by the Baire category theorem, every countable compact homogeneous space is finite. Hence we reach a contradiction. Thus $K$ is not discrete. Since $K$ is a closed subgroup of the Lie group $H$, it is itself a Lie group, hence $\operatorname{Lie}(K)\neq0$. Since $G_x$ is discrete, $\operatorname{Lie}(G_x) = 0$, and $K\subseteq H_Y\cap(G_x\times G)$, we get 
\[
0\neq\operatorname{Lie}(K) \subseteq \mathfrak h_Y\cap(0\oplus\mathfrak g),
\]
where $\mathfrak h_Y:=\operatorname{Lie}(H_Y)$. In other words, there exists $Z\in\mathfrak g\setminus\{0\}$, such that $(0,Z)\in \mathfrak h_Y\cap(0\oplus\mathfrak g)$. Thus the set 
\[
\mathfrak i:= \{Z\in\mathfrak g:(0,Z)\in\mathfrak h_Y\}
\]
is nonzero. Since $Y$ is a finite-volume homogeneous subspace and is $\Sigma$-invariant, $\Sigma$-ergodic, by definition $\Sigma$ is a sub-semi-group of the stabilizer $H_Y$ and acts ergodically on $(Y,\nu_Y)$. In particular, $\Sigma\subseteq H_Y$. Hence it follows that, $(\gamma,\gamma)\in H_Y$ for every $\gamma\in\Gamma$. Thus the adjoint action of $(\gamma,\gamma)$ preserves the Lie algebra $\mathfrak h_Y$. In other words, if $Z\in\mathfrak i$, then $(0,Z)\in\mathfrak h_Y$ and
\[
\Ad_H(\gamma,\gamma)(0,Z) = (0,\Ad(\gamma)Z) \in \mathfrak h_Y.
\]
Hence $\Ad(\gamma)Z\in\mathfrak i$ for every $\gamma\in\Gamma$, and therefore $\mathfrak i$ is invariant under $\Ad(\Gamma)$. By the $\Ad$-Zariski density of $\Gamma$, $\mathfrak i$ is invariant under $\Ad(G)$. Differentiating the adjoint action gives
\[
[\mathfrak g,\mathfrak i]\subseteq\mathfrak i.
\]
Thus $\mathfrak i$ is an ideal of $\mathfrak g$. Since $\mathfrak g$ is simple and $\mathfrak i\neq0$, we get $\mathfrak i=\mathfrak g$. Hence by definition of $\mathfrak i$, we have 
\begin{equation}\label{eq:0_g_h_y}
    0\oplus\mathfrak g\subseteq\mathfrak h_Y.
\end{equation}
Since $G$ is a connected Lie group, \eqref{eq:0_g_h_y} implies 
\[
\{e\}\times G\subseteq H_Y.
\]
Since the homogeneous measure $\nu_Y$ is invariant under the action of $H_Y$, it follows that $\nu_Y$ is invariant under the action of $\{e\}\times G$. We define a finite complex measure $\rho$ on $X$ as follows. Let
\[
\int_X\psi(v)\,d\rho(v) := \int_{X\times X}f(u)\psi(v)\,d\nu_Y(u,v),
\] for every $\psi\in C(X)$. This defines $\rho$ by the Riesz representation theorem. The invariance of $\nu_Y$ under $\{e\}\times G$ implies that $\rho$ is $G$-invariant. Since every finite $G$-invariant complex measure on $G/\Lambda$ is a scalar multiple of $m$, it follows that $\rho=c\,m$ for some $c\in\C$. Taking $\psi=\overline f$, we obtain
\[
\int_Y\Phi_f\,d\nu_Y = \int_{X\times X}f(u)\overline{f(v)}\,d\nu_Y(u,v) = \int_X\overline{f(v)}\,d\rho(v) = c\int_X\overline f\,dm = 0. 
\]
This contradicts \eqref{eq:_Phi_contra}, hence $E_{f,\eta}(x)$ is finite. Thus the sequence $(\mu_N)_{N\geq1}$ satisfies the hypothesis of \Cref{thm:main}. Hence \Cref{thm:main} implies that $\Gamma\curvearrowright X$ has the Glasner property and is IC-rigid.
\end{proof}

Note that the \Cref{cor:homogeneous} follows immediately from \Cref{thm:simple}. 

\begin{proof}[Proof of \Cref{cor:homogeneous}] 
For $G=\mathrm{SL}_d(\R)$ or $G=\mathrm{PSL}_d(\R)$, the group $G$ is a connected noncompact Lie group and
\[
\operatorname{Lie}(G)=\mathfrak{sl}_d(\R),
\]
which is simple for $d\geq2$. Since $\Gamma$ is Zariski dense in $G$, and the adjoint representation is an algebraic homomorphism  it follows that $\Ad(\Gamma)$ is Zariski dense in $\Ad(G)$. Thus $\Gamma$ is $\Ad$-Zariski dense in $G$. Therefore all the hypotheses of \Cref{thm:simple} are satisfied. Hence the action $\Gamma\curvearrowright G/\Lambda$ has the Glasner property and is IC-rigid.
\end{proof}

\section{Glasner Property and IC-rigidity of Torus actions}\label{sec:Torus}

In this section we prove \Cref{thm:Fourier} which is a Fourier-analytic condition that implies IC-rigidity and the Glasner property for semi-group actions on $d$-dimensional torus. Using this criterion we prove \Cref{thm:polynomial-matrix}, i.e. the matrix-valued polynomial actions on the $d$-dimensional torus are IC-rigid. The Glasner property for these actions was previously established by Bulinski and Fish \cite{Bulinski_Fish_2023_a}, who in fact proved the stronger uniform Glasner property. 

To do so we need an elementary result showing that small Fourier coefficients for a measure on the torus force dense support (see \Cref{prop:Fourier-density}). Note that this result can be derived from \cref{prop:density-detection} but we provide a simple standalone proof here. Recall that for a  Borel probability measure $\mu$ on $\T^d$, and $m \in \Z^d$, we denote the $m$-th Fourier coefficient of $\mu$ by
\[
\widehat\mu(m) :=\int_{\T^d}e\bigl(\langle m,x\rangle\bigr)\,d\mu(x).
\]
\begin{prop}\label{prop:Fourier-density}
For every $\varepsilon>0$, there exist a constant $c_\varepsilon>0$ and a finite nonempty set $\Lambda_\varepsilon\subseteq\Z^d\setminus\{0\}$ such that the following holds. If $\mu$ is a Borel probability measure on $\T^d$ satisfying
\[
|\widehat\mu(m)|<c_\varepsilon
\qquad\text{for every }m\in\Lambda_\varepsilon,
\]
then $\operatorname{supp}\mu$ is $\varepsilon$-dense in $\T^d$.
\end{prop}

\begin{proof} Let $\lambda$ denote normalized Haar measure on $\T^d$. Choose $y_1,\ldots,y_s\in\T^d$ such that the collection of open balls $\{U_j:=B(y_j,\varepsilon/2): 1 \leq j \leq s\}$, cover $\T^d$. For each $1 \leq j \leq s$, choose a non-negative continuous function $f_j$ such that $\operatorname{supp}f_j\subseteq U_j$ and $\int_{\T^d}f_j\,d\lambda=1$. Such a function can be obtained by taking a non-negative bump function supported in a smaller ball around $y_j$ and dividing by its integral. Using the fact that trigonometric polynomials are uniformly dense in $C(\T^d)$, for every $1 \leq j \leq s$, we can choose a trigonometric polynomial $P_j$ such that 
\[\|P_j-f_j\|_\infty<\frac{1}{8}.\] Note that there exists a common finite set $\Lambda_\varepsilon \subseteq \Z^d \setminus \{0\}$ such that for every $1 \leq j \leq s$ we can write 
\[
P_j(x)=a_{j,0} +\sum_{m\in\Lambda_\varepsilon} a_{j,m}e\bigl(\langle m,x\rangle\bigr)
\]
for $x \in \T^d$. We denote by
\[
C:=\max_{1\leq j\leq s} \sum_{m\in\Lambda_\varepsilon}|a_{j,m}|, \qquad\text{and }\qquad c_\varepsilon:=\frac{1}{4(1+C)}.
\] Since $\lambda$ is a Haar measure, for $m\in\Z^d \setminus \{0\}$, $\int_{\T^d} e\bigl(\langle m,x\rangle\bigr)\,d\lambda(x)=0$. Thus for every $1 \leq j \leq s$,
\[
\int_{\T^d}P_j\,d\mu-\int_{\T^d}P_j\,d\lambda
=\sum_{m\in\Lambda_\varepsilon}a_{j,m} \, \widehat\mu(m).
\] Using the fact that $|\widehat\mu(m)|<c_\varepsilon$ it follows that for every $1 \leq j \leq s$,
\[
\left|\int_{\T^d}f_j\,d\mu-1\right| = \left|\int_{\T^d}f_j\,d\mu -\int_{\T^d}f_j\,d\lambda\right| \leq 2\|f_j-P_j\|_\infty + \sum_{m\in\Lambda_\varepsilon}|a_{j,m}|\,|\widehat\mu(m)| < \frac14+c_\varepsilon C < \frac12.
\]
Thus $\int_{\T^d}f_j\,d\mu>1/2$ for every $j \in \{1,\ldots,s\}$. Since $f_j$ vanishes outside $U_j$ and $\int_{\T^d} f_j\,d\mu>0,$
it follows that $\operatorname{supp}\mu\cap U_j\neq\varnothing$ for every $j \in \{1,\ldots,s\}$. Since $U_j$ are balls of radius $\varepsilon/2$ that cover $\T^d$, it follows that $\operatorname{supp}\mu$ is $\varepsilon$-dense in $\T^d$. 
\end{proof}

Now we are ready to prove \Cref{thm:Fourier}.

\begin{proof}[Proof of \Cref{thm:Fourier}]
Let $\lambda$ denote normalized Haar measure on $\T^d$. Fix $\varepsilon>0$ and let $c_\varepsilon>0$ and $\Lambda_\varepsilon\subseteq\Z^d\setminus\{0\}$ be as in \Cref{prop:Fourier-density}. We set $q:=|\Lambda_\varepsilon|$ and enumerate $\Lambda_\varepsilon=\{m_1,\ldots,m_q\}$. For every $1\leq j\leq q$ we define $f_j(x):=e\bigl(\langle m_j,x\rangle\bigr)$. Since $m_j\neq0$, it follows that $\int_{\T^d}f_j\,d\lambda=0$, for every $1\leq j\leq q$. In other words $f_1,\ldots,f_q\in C(\T^d;\C)$, is a collection of zero-mean functions. Note that for every Borel probability measure $\nu$ on $\T^d$,
\[
\int_{\T^d}f_j\,d\nu =\int_{\T^d}e\bigl(\langle m_j,x\rangle\bigr)\,d\nu(x) = \widehat\nu(m_j).
\]
Hence, by \Cref{prop:Fourier-density}, if
\[
\left|\int_{\T^d}f_j\,d\nu\right|<c_\varepsilon
\qquad\text{for every }1\leq j\leq q,
\]
then $\operatorname{supp}\nu$ is $\varepsilon$-dense in $\T^d$. Thus the functions $f_1,\ldots,f_q$, with the common constant $c_\varepsilon$, have precisely the $\varepsilon$-density property used in the proof of \Cref{thm:main}. Recall that in the proof of \Cref{thm:main}, for each $\varepsilon>0$, we used \Cref{prop:density-detection} to obtain a finite collection of zero-mean functions and positive constants which detect $\varepsilon$-density.  Here \Cref{prop:Fourier-density} provides such a collection, namely $f_1,\ldots,f_q$, with common constant $c_\varepsilon$. The remainder of the proof of \Cref{thm:main} uses the pair-correlation condition only for this finite collection. We now verify the pair-correlation condition for this finite collection of functions using the hypothesis of \Cref{thm:Fourier}. For $1\leq j\leq q$ and $x,y\in\T^d$ we have
\[
K_{N,f_j}(x,y) = \int_G e\bigl(\langle m_j,gx\rangle\bigr)\overline{e\bigl(\langle m_j,gy\rangle\bigr)}\,d\mu_N(g) = \int_G e\bigl(\langle m_j,g(x-y)\rangle\bigr)\,d\mu_N(g).
\]
Fix $x\in\T^d$ and $\eta>0$, and set
\[
E_{m_j,\eta} := \left\{ z\in\T^d: \limsup_{N\to\infty} \left| \int_G e\bigl(\langle m_j,gz\rangle\bigr)\,d\mu_N(g) \right| \geq\eta \right\}.
\]
By the hypothesis of \Cref{thm:Fourier}, $E_{m_j,\eta}$ is finite. If
\[
\limsup_{N\to\infty}|K_{N,f_j}(x,y)|>\eta,
\]
then $x-y\in E_{m_j,\eta}$. Hence it follows that 
\[
R_{j,\eta}(x) : = \left\{ y\in\T^d\setminus\{x\}: \limsup_{N\to\infty}|K_{N,f_j}(x,y)|>\eta \right\} \subseteq x-E_{m_j,\eta}. 
\]
Thus $R_{j,\eta}(x)$ is finite for every $1\leq j\leq q$. Hence the proof of \Cref{thm:main}, applied with the finite collection of functions above, gives that the action $G\curvearrowright\T^d$ has the Glasner property and is IC-rigid.
\end{proof}

\subsection{Matrix-valued polynomial actions are IC-rigid}

\begin{proof}[Proof of \Cref{thm:polynomial-matrix}]
Let $A(x)\in M_d(\Z[x])$ be as in the statement of the theorem. We consider the action of the semi-group
\[
G_A = \left\{ A(n_1)A(n_2)\cdots A(n_\ell): \ell\geq 1,\ n_1,\ldots,n_\ell\in\Z,\ \det A(n_i)\neq 0 \text{ for all }1\leq i\leq \ell \right\}
\]
on $\T^d$ by endomorphisms. Since $\det A(x)$ is a nonzero polynomial, there are only finitely many integers $n$ for which $\det A(n)=0$. Thus, for all sufficiently large $N$, let
\[
I_N:=\{n\in\Z:|n|\leq N,\ \det A(n)\neq0\}
\]
and define a probability measure on $G_A$ by
\[
\mu_N:=\frac{1}{|I_N|} \sum_{n\in I_N}\delta_{A(n)}.
\] For the finitely many remaining $N$, define $\mu_N$ arbitrarily. Since only finitely many integers are excluded, it follows that $\frac{|I_N|}{2N+1}\rightarrow1$. For $z\in\T^d$, the average orbit probability measure $\sigma_{N,z}$
(corresponding to $\mu_N$) is given by
\[
\sigma_{N,z}
:=
\int_{G_A}\delta_{gz}\,d\mu_N(g)
=
\frac1{|I_N|}
\sum_{n\in I_N}\delta_{A(n)z}.
\]
For $m\in\Z^d$, the $m$-th Fourier coefficient of $\sigma_{N,z}$ is given by
\[
\widehat{\sigma_{N,z}}(m) = \int_{\T^d}e(\langle m,x\rangle)\,d\sigma_{N,z}(x) = \int_{G_A}e(\langle m,gz\rangle)\,d\mu_N(g) = \frac1{|I_N|} \sum_{n\in I_N} e\bigl(\langle m,A(n)z\rangle\bigr).
\]
Fix $m\in\Z^d\setminus\{0\}$ and $\delta>0$, and set
\[
E_{m,\delta} := \left\{ z\in\T^d: \limsup_{N\to\infty} \left| \frac1{|I_N|} \sum_{n\in I_N} e\bigl(\langle m,A(n)z\rangle\bigr) \right| \geq\delta \right\}.
\]
To prove that the action is IC-rigid, by \Cref{thm:Fourier} (b), it is enough to show that $E_{m,\delta}$ is finite. 

We first consider a non-torsion point $z\in\T^d$. For
$m\in\Z^d$ let $m^t$ denote the corresponding row vector. We write
\[
m^t(A(x)-A(0))=(P_1(x),\ldots,P_d(x)),
\]
where $P_1,\ldots,P_d\in\Z[x]$. The polynomials $P_1,\ldots,P_d$ are linearly independent over $\Z$. Indeed, otherwise there would exist $w\in\Z^d\setminus\{0\}$ such that
\[
m\cdot(A(x)-A(0))w=0,
\]
contrary to the hypothesis of the theorem. Since $z=(z_1,\ldots,z_d)$ is non-torsion, at least one $z_i$ is irrational. We claim that the polynomial $\sum_{i=1}^d z_iP_i(x)$ has an irrational coefficient. To see this note that, since $P_1,\ldots,P_d$ are linearly independent over $\Z$ they are also linearly independent over $\Q$. Let $D$ be the maximum of the degrees of the entries of $A(x)$. Since $P_i(0)=0$ and $\deg P_i\leq D$ for every $i$, we may write $P_i(x)=\sum_{j=1}^D b_{ji}x^j$, where $b_{ji}\in\Z$. 

Let $B\in M_{D\times d}(\Z)$ be the matrix whose $i$-th column is the coefficient vector of $P_i$. Since $P_1,\ldots,P_d$ are linearly independent over $\Q$, the matrix $B$ has full column rank over $\Q$. Assume by contradiction that polynomial $\sum_{i=1}^d z_iP_i(x)$ has only rational coefficients then $Bz\in\Q^D$. Since $B$ has full column rank over $\Q$, this would imply that $z\in\Q^d$, contradicting the assumption that $z$ is non-torsion. Hence $\sum_{i=1}^d z_iP_i(x)$ has an irrational coefficient. Since $P_i(0)=0$ for every $i$, this coefficient is necessarily non-constant. Note that
\[
\langle m,A(x)z\rangle = m^tA(0)z + m^t(A(x)-A(0))z = \langle m,A(0)z\rangle + \sum_{i=1}^d z_iP_i(x),
\]
and the first term is constant in $x$. Hence $\langle m,A(x)z\rangle$ has at least one irrational non-constant coefficient. By Weyl's equidistribution theorem, we get
\[
\lim_{N\to\infty} \frac{1}{2N+1} \sum_{|n|\leq N} e(\langle m,A(n)z\rangle) =0.
\]
Removing the finitely many $n$ for which $\det A(n)=0$ does not affect the limit. Therefore
\[
\lim_{N\to\infty} \frac{1}{|I_N|} \sum_{n\in I_N} e(\langle m,A(n)z\rangle) =0.
\]
Thus every point of $E_{m,\delta}$ is a torsion point. For a torsion point $z\in\T^d$, let
\[
z=\frac{a}{q}+\Z^d,\qquad a=(a_1,\ldots,a_d)\in\Z^d,\qquad \gcd(a_1,\ldots,a_d,q)=1
\]
be its reduced form. Since $A(x)$ has integer polynomial entries and $m,a\in\Z^d$, the quantity $m\cdot A(x)a$ is a polynomial in $\Z[x]$ of degree at most $D$. Write
\[
m\cdot A(x)a=b_0+b_1x+\cdots+b_Dx^D,
\]
where $b_0,\ldots,b_D\in\Z$. Using the fact that $m\cdot A(0)a=b_0$, we get
\begin{equation}\label{eq:m^tA}
\bigl(m^tA(x)-m^tA(0)\bigr)a=b_1x+\cdots+b_Dx^D.
\end{equation}
On the other hand, since $m^t(A(x)-A(0))=(P_1(x),\ldots,P_d(x))$ the left-hand side in \eqref{eq:m^tA} is equal to $\sum_{i=1}^d a_iP_i(x)$. By the definition of $B$, the coefficient vector of this polynomial is $Ba$. Therefore
\[
Ba=(b_1,\ldots,b_D)^t.
\]
Since $B$ has full column rank over $\Q$, there exists a matrix $L\in M_{d\times D}(\Q)$ such that $LB=I_d$. 
Clearing denominators, there exist $C\in M_{d\times D}(\Z)$ and an integer $Q_m\geq 1$, depending only on $A(x)$ and $m$, such that $CB=Q_mI_d$. 
Let $r:=\gcd(b_1,\ldots,b_D,q)$. Since $r$ divides every entry of $Ba$ and also divides $q$, it follows from $CBa=Q_ma$
that $r\mid \gcd(Q_ma_1,\ldots,Q_ma_d,q)$. Since $\gcd(a_1,\ldots,a_d,q)=1$, we have $\gcd(Q_ma_1,\ldots,Q_ma_d,q)\mid Q_m$. Hence $r\leq Q_m$. We set
\[
q_0:=\frac{q}{r}, \qquad b_j':=\frac{b_j}{r}, \quad 1\leq j\leq D.
\]
Then
\begin{equation}\label{eq:GCD_b'_q_0}
    \gcd(b_1',\ldots,b_D',q_0)=1.
\end{equation}
Since $m\cdot A(x)a$ has integer coefficients, $m\cdot A(n+q)a-m\cdot A(n)a$ is divisible by $q$. Hence the sequence
\[
e(\langle m,A(n)z\rangle) = e\left(\frac{m\cdot A(n)a}{q}\right)
\]
is periodic in $n$ with period $q$. Since only finitely many terms are removed from the symmetric averages, we obtain
\[
\lim_{N\to\infty} \frac{1}{|I_N|} \sum_{n\in I_N} e(\langle m,A(n)z\rangle) = \frac{1}{q} \sum_{n=1}^{q} e\left(\frac{m\cdot A(n)a}{q}\right).
\]
Since $b_j=rb_j'$ for $1\leq j\leq D$ and $q=rq_0$, we get
\[
\frac{m\cdot A(n)a}{q} = \frac{b_0}{q} + \frac{b_1'n+\cdots+b_D'n^D}{q_0}.
\]
Hence
\[
e\left(\frac{m\cdot A(n)a}{q}\right) = e\left(\frac{b_0}{q}\right) e\left( \frac{b_1'n+\cdots+b_D'n^D}{q_0}
\right).
\]
The first factor is independent of $n$ and has absolute value one. Note that the numerator of the second factor is a polynomial with integer coefficients, hence the second factor is periodic in $n$ with period $q_0$. Since $q=rq_0$, the sum over $1\leq n\leq q$ consists of $r$ identical blocks of length $q_0$. Therefore
\begin{equation}\label{eq:m_A_n_a}
    \left| \frac{1}{q} \sum_{n=1}^{q} e\left(\frac{m\cdot A(n)a}{q}\right) \right| = \left| \frac{1}{q_0} \sum_{n=1}^{q_0} e\left( \frac{b_1'n+\cdots+b_D'n^D}{q_0} \right) \right|.
\end{equation}
Recall that by Hua’s estimate for exponential sums \cite[Chapter I, Theorem 1]{Hua_1965}, if
\[
f(x)=A_kx^k+\cdots+A_1x+A_0 \in \Z[x]
\]
and $\gcd(A_k,\ldots,A_1,q)=1$ then for every $\eta>0$, there exists $C_{k,\eta} > 0$ such that
\[
\left| \sum_{n=1}^{q} e\left(\frac{f(n)}{q}\right) \right| \leq C_{k,\eta}q^{1-\frac1k+\eta}.
\]
We apply this with $f(x)=(b_1'x+\cdots+b_D'x^D)$ and $q=q_0$. Let $k\leq D$ be the degree of $f$. Using \eqref{eq:GCD_b'_q_0} and dividing both sides by $q_0$, we get
\[
\left| \frac{1}{q_0} \sum_{n=1}^{q_0} e\left( \frac{b_1'n+\cdots+b_D'n^D}{q_0} \right) \right| \leq C_{k,\eta}q_0^{-\frac1k+\eta}.
\]
Let $\eta=1/(2D)$. Since $k\leq D$, we have $-\frac{1}{k}+\frac{1}{2D}
\leq -\frac{1}{2D}.$ Absorbing the finitely many possible values of $k$ into a constant $C_D$, and applying \eqref{eq:m_A_n_a} we get
\[
\left| \frac{1}{q} \sum_{n=1}^{q} e\left(\frac{m\cdot A(n)a}{q}\right) \right| \leq C_Dq_0^{-\frac{1}{2D}}.
\]
Since $q_0=q/r$ and $r\leq Q_m$,
\[
\left| \frac{1}{q} \sum_{n=1}^{q} e\left(\frac{m\cdot A(n)a}{q}\right) \right| \leq C_D \left(\frac{Q_m}{q}\right)^{\frac{1}{2D}}.
\]
Thus, if $z\in E_{m,\delta}$, then $\delta \leq C_D \left(\frac{Q_m}{q}\right)^{\frac{1}{2D}}$ and consequently
\[
q \leq Q_m\left(\frac{C_D}{\delta}\right)^{2D}.
\]
Hence the order of every point in $E_{m,\delta}$ is bounded in terms of $m$, $\delta$, and $A(x)$. There are only finitely many torsion points of $\T^d$ of bounded order, and therefore $E_{m,\delta}$ is finite. The hypothesis of \Cref{thm:Fourier} is satisfied. Hence, by \Cref{thm:Fourier} (b), the action $G_A\curvearrowright\T^d$ is IC-rigid.
\end{proof}

\bibliographystyle{alpha}
\bibliography{references}

\end{document}